\documentclass[11pt]{article}
\usepackage[utf8]{inputenc}
\usepackage[english]{babel}
\usepackage{times}
\usepackage{paralist}
\usepackage{amsmath,amstext,amssymb,amsthm, amsfonts}
\usepackage{graphicx}

\newtheorem{lemma}{Lemma}[section]
\newtheorem{theorem}[lemma]{Theorem}
\newtheorem{corollary}[lemma]{Corollary}
\newtheorem{definition}[lemma]{Definition}

\newtheorem{example}[lemma]{Example}

\newtheorem{Assumption}[lemma]{Assumption}

\newcommand{\ilim} {\mathop{\rm lim\,inf\,}}

\def\U{\mathbb{U}}

\def\S{\mathbb{S}}
\def\K{\mathbb{K}}
\def\M{\mathbb{M}}
\def\Y{\mathbb{Y}}
\def\h{\mathbf{I}}
\def\X{\mathbb{X}}
\def\E{\mathbb{E}}
\def\A{\mathbb{A}}
\def\H{\mathbb{H}}
\def\R{\mathbb{R}}
\def\P{\mathbb{P}}
\def\F{\mathbb{F}}

\def\C{\mathbb{C}}
\def\W{\mathbb{W}}
\def\lll{\mathbb{L}}
\def\c{\bar{c}}
\def\B{\mathcal{B}}
\def\oo{\mathcal{O}}

\def\fff{\mathtt{F}}

\def\Psii{{\Xi_{\int}}}

\title{A Class of Solvable Markov Decision Models with Incomplete Information}
\begin{document}

\maketitle

\begin{center}
Eugene~A.~Feinberg \footnote{Department of Applied Mathematics and
Statistics,
 Stony Brook University,
Stony Brook, NY 11794-3600, USA, eugene.feinberg@sunysb.edu},\
Pavlo~O.~Kasyanov\footnote{Institute for Applied System Analysis,
National Technical University of Ukraine ``Kyiv Polytechnic
Institute'', Peremogy ave., 37, build, 35, 03056, Kyiv, Ukraine,\
kasyanov@i.ua.},\ and Michael~Z.~Zgurovsky\footnote{National Technical University of Ukraine
``Kyiv Polytechnic Institute'', Peremogy ave., 37, build, 1, 03056,
Kyiv, Ukraine,\
mzz@kpi.ua
}\\

\bigskip
\end{center}

\begin{abstract}
This paper investigates natural conditions for the existence of optimal policies for a  Markov decision process with incomplete information (MDPII) and with expected total costs.  The MDPII is the classic model of a controlled stochastic process with incomplete state observations which is more general than  Partially Observable Markov Decision Processes (POMDPs). For MDPIIs we introduce the notion of a semi-uniform Feller transition probability, which is stronger than the notion of a weakly continuous transition probability.  We show that an MDPII has a semi-uniform Feller transition probability if and only if the corresponding belief MDP also has a semi-uniform Feller transition probability.  This fact has several corollaries. In particular, it provides new and implies all known sufficient conditions for the existence of optimal policies for POMDPs with expected total costs.
\end{abstract}

\section{Introduction} In many control problems the state of a controlled system is not known, and decision makers know only some information about the state. This takes place in many applications including signal processing, robotics, artificial intelligence, and medicine.  Except lucky exceptions, and Kalman's filtering is among them, problems with incomplete information are known to be difficult \cite{PT}.   The general approach to solving such problems was identified long time ago in \cite{Ao, As, Dy, ShP}, and it is based on constructing a controlled system whose states are posterior state distributions for the original system.  These posterior distributions are often called belief probabilities or belief states.  Finding an optimal policy for a problem with incomplete state observation consists of two steps: (i) finding an optimal policy for the problem with belief states, and (ii) deriving  an optimal policy for the original problem  from this policy. This approach was introduced in \cite{Ao, As, Dy, ShP} for problems with finite state, observation, and action sets, and it holds for problems with Borel state, observation, and action sets \cite{Rh, Yu}.  If there is no optimal policy for the problem with belief space, then there is no optimal policy for the original problem.

  This paper deals with optimization of  expected total discounted costs for models with discrete time.  We describe a large class of problems, for which optimal policies exist, satisfy optimality equations, and can be found by value iterations. For a particular model of Partially Observable Markov Decision Processes (POMDPs) the related studies are \cite{FKZ, Saldi}.

There are several models of controlled systems with incomplete state observations in the literature. Here we mostly consider a contemporary version of the original model introduced in  \cite{Ao, As, Dy, ShP} and called a Markov Decision Process with Incomplete Information (MDPII).  In this model the transitions are defined by transition probabilities $P(dw_{t+1},dy_{t+1}|w_t,y_t,a_t),$  where vectors $(w_t,y_t)$ represent states of the system at times $t=0,1,\ldots,$  $w_t$ and $y_t$ are unobservable and observable components of the state $(w_t,y_t)$, and $a_t$  are actions. In more contemporary studies the research focus switched to POMDPs. As was observed in \cite{Plat}, there are two different POMDP models in the literature, which we call ${\rm POMDP}_1$ and ${\rm POMDP}_2.$  Platzman~\cite{Plat} introduced a ``plant'' model, which we call Platzman's model, and this model is more general than ${\rm POMDP}_1$ and ${\rm POMDP}_2$.

 In order to explain that
 Platzman's plant is a particular case of an MDPII, let us change the notation $w_t$ to $x_t$ in this paragraph, that is $x_t:=w_t,$ and this change takes place only in this paragraph.  Then the transition probability for an MDP is
 $P(dx_{t+1},dy_{t+1}|x_t,y_t,a_t),$ and $x_t$ and $y_t$ can be interpreted as a hidden state  and  an observation.  Platzman's model is a particular case of an MDPII when the transition probability does not depend on observations.
 In other words, the transition probability in Platzman's model is $P(dx_{t+1},dy_{t+1}|x_t,a_t).$   ${\rm POMDP}_i$ are Platzman’s models whose transition probabilities have special structural properties. These properties are  $P(dx_{t+1},dy_{t+1}|x_t,a_t)=P_1(dx_{t+1}| x_t,a_t)$ $Q_1(dy_{t+1}|x_t,a_t)$ for ${\rm POMDP}_1$ and $P(dx_{t+1},dy_{t+1}|x_t,a_t)= Q_2(dy_{t+1}|a_t,x_{t+1}) P_2(dx_{t+1}| x_t,a_t)$ for ${\rm POMDP}_2.$

This paper introduces the class of MDPIIs with semi-uniform Feller transition probabilities.  We show that, if an MDPII has a transition probability from this class, then the transition probability of the belief-MDP also belongs to this class.  Under mild conditions on cost functions, there are optimal policies for MDPs with semi-uniform Feller transition probabilities.  This paper provides several sufficient conditions for the existence of optimal policies, validity of optimality equations, and convergence of value iterations.  This general theory implies the following sufficient conditions for the existence of optimal policies for POMDPs: (i) $P_i$ is weakly continuous and $Q_i$  is continuous in total variation for  a $ {\rm POMDP}_i,$ $i=1,2$  (for $i=2$ this result was established in \cite{FKZ}); (ii) $P_2$ is continuous in total variation and $Q_2$ is continuous in total variation in the control parameter; sufficiency of continuity of $P_2$ in total variation was established in \cite{Saldi} for uncontrolled observation kernels, that is, $ Q_2(y_{t+1}|a_t,x_{t+1})= Q_2(y_{t+1}|x_{t+1})$.

The results are presented without proofs.  The proofs can be found in \cite{FKZ2021a, FKZ2021b}.

\section{Model Description}\label{sec:MoDe}

For a metric space $\S=(\S,\rho_\S),$ where $\rho_\S$ is a metric, let $\tau(\S)$ be the topology of $\S$ (the family of all open subsets of $\S$), and let ${\mathcal B}(\S)$ be its Borel
$\sigma$-field, that is, the $\sigma$-field generated by all open subsets of the
metric space $\S$. For $s\in \S$ and $\delta>0$ denote by
$B_\delta(s)$ and $\bar{B}_\delta(s)$ respectively the open and
closed balls in the
metric space $\S$ of  radius $\delta$ with  center $s$ and by $S_\delta(s)$ the sphere in
$\S$ of  radius $\delta$ with  center $s$. Note that $S_\delta(s)=\bar{B}_\delta(s)\setminus B_\delta(s).$
For a subset $S$ of $\S$ let $\bar{S}$ denote the \textit{closure of} $S$ and $S^o$ the \textit{interior of} $S.$ Then $S^o\subset S\subset \bar{S}.$ $S^o$ is open and $\bar{S}$ is closed. $\partial S:=\bar{S}\setminus S^o$ denotes the \textit{boundary of} $S.$ In particular, $\partial B_\delta(s)=S_\delta(s).$
We denote by $\P(\S)$ the \textit{set of probability
measures} on $(\S,{\mathcal B}(\S)).$ A sequence of probability
measures $\{\mu^{(n)}\}_{n=1,2,\ldots}$ from $\P(\S)$
\textit{converges weakly} to $\mu\in\P(\S)$ if for any
bounded continuous function $f$ on $\S$
\[\int_\S f(s)\mu^{(n)}(ds)\to \int_\S f(s)\mu(ds) \qquad {\rm as \quad
}n\to\infty.
\]
A sequence of probability measures $\{\mu^{(n)}\}_{n=1,2,\ldots}$ from $\P(\S)$ \textit{converges in  total  variation} to $\mu\in\P(\S)$ if
\begin{equation}\label{eq:Kara1}
\begin{aligned}
\sup_{C\in \B(\S)}|\mu^{(n)}(C)-\mu(C)|\to 0\ {\rm as} \ n\to \infty;
\end{aligned}
\end{equation}
see~\cite{UFL} for properties of these types of convergence of probability measures.
Note that $\P(\S)$ is a separable metric space with respect to the topology of weak convergence of probability measures, when $\S$ is
a separable metric space; \cite[Chapter~II]{Part}. Moreover,
according to Bogachev \cite[Theorem~8.3.2]{bogachev}, if the metric space $\S$ is separable, then the topology of weak convergence of probability measures on $(\S,\B(\S))$ coincides with the topology generated by the \textit{Kantorovich-Rubinshtein metric}
\begin{equation}\label{eq:KantorRubMetr}
\begin{aligned}
\rho_{\P(\S)}(\mu,\nu):=\sup  \left\{\int_\S\right. &f(s)\mu(ds)-\int_\S f(s)\nu(ds) \ \Big{|} \\ &\left.f\in{\rm Lip}_1(\S),\ \sup_{s\in\S}|f(s)|\le 1 \right\},
\end{aligned}
\end{equation}
$\mu,\nu\in\P(\S),$ where ${\rm Lip}_1(\S):=\{f:\S\mapsto\mathbb{R}, \ |f(s_1)-f(s_2)|\le \rho_\S(s_1,s_2),\ \forall s_1,s_2\in\S\}.$

For a Borel subset $S$ of a metric space $(\S,\rho_\S)$, where $\rho_\S$
is a metric, we always consider the
 metric space $(S,\rho_S),$ where $\rho_S:=\rho_\S\big|_{S\times S}.$  A subset $B$ of $S$ is called open
(closed) in $S$ if $B$ is open (closed respectively) in $(S,\rho)$. Of course, if $S=\S$, we omit
``in $\S$''. Observe that, in general, an open (closed) set
in $S$ may not be open (closed respectively). For $S\in\B(\S)$ we denote by $\B(S)$ the Borel $\sigma$-field on
$(S,\rho_S).$  Observe that $\B(S)=\{S\cap B:B\in\B(\S)\}.$

For metric spaces $\S_1$ and $\S_2$, a (Borel-measurable) \textit{stochastic
kernel} $\Psi(ds_1|s_2)$ on $\S_1$ given $\S_2$ is a mapping $\Psi(\,\cdot\,|\,\cdot\,):\B(\S_1)\times \S_2\mapsto [0,1]$, such that $\Psi(\,\cdot\,|s_2)$ is a
probability measure on $\S_1$ for any $s_2\in \S_2$, and $\Psi(B|\,\cdot\,)$ is a Borel-measurable function on $\S_2$ for any Borel set $B\in\B(\S_1)$.  Another name for a stochastic kernel is transition probability. A
stochastic kernel $\Psi(ds_1|s_2)$ on $\S_1$ given $\S_2$ defines a Borel measurable mapping $s_2\mapsto \Psi(\,\cdot\,|s_2)$ of $\S_2$ to the metric space
$\P(\S_1)$ endowed with the topology of weak convergence.
A stochastic kernel
$\Psi(ds_1|s_2)$ on $\S_1$ given $\S_2$ is called
\textit{weakly continuous (continuous in
total variation)}, if $\Psi(\,\cdot\,|s^{(n)})$ converges weakly (in
total  variation) to $\Psi(\,\cdot\,|s)$ whenever $s^{(n)}$ converges to $s$
in $\S_2$. For one-point sets $\{s_1\}\subset \S_1,$ we
sometimes write $\Psi(s_1|s_2)$ instead of $\Psi(\{s_1\}|s_2)$. Sometimes a weakly continuous stochastic kernel is called Feller, and a stochastic kernel continuous in total variation is called uniformly Feller \cite{Papa}.

Let $\S_1,\S_2,$ and $\S_3$ be Borel subsets of Polish spaces (a Polish space is a complete
separable metric space), and $\Psi$ on $\S_1\times\S_2$ given $\S_3$ be a stochastic kernel. For each $A\in\B(\S_1),$
$B\in\B(\S_2),$ and $s_3\in\S_3,$ let:
\begin{equation}\label{eq:marg_new}
\Psi(A,B|s_3):=\Psi(A\times B|s_3).
\end{equation}
In particular, we consider \textit{marginal} stochastic kernels
$\Psi(\S_1,\,\cdot\,|\,\cdot\,)$ on $\S_2$ given $\S_3$ and $\Psi(\,\cdot\,,\S_2|\,\cdot\,)$ on $\S_1$ given $\S_3.$

A \textit{Markov decision process with incomplete
information (MDPII)} is specified by a tuple
$(\W\times\Y,\A,P,c),$ where

(i) $\W\times\Y$ is the \textit{state space}, where $\W$ and $\Y$ are Borel subsets of Polish spaces, and for $(w,y)\in\W\times\Y,$ the unobservable component
    of the state $(w,y)$ is $w$ and the observable component is $y;$

(ii) $\A$ is the \textit{action space}, which is assumed to be a Borel subset of a Polish space;

(iii) $P$ is a stochastic kernel on $\W\times\Y$ given
$\W\times\Y\times\A,$ which determines the distribution $P(\,\cdot\,|w,y,a)$ on $\W\times\Y$ of the new state, if $(w,y)\in\W\times\Y$ is the current state, and if $a\in A(y)$ is the current action, and it is assumed that the stochastic  kernel $P$ on $\W\times\Y$ given $\W\times\Y\times\A$
is weakly continuous in $(w,y,a)\in \W\times\Y\times\A;$

(iv) $P_0(\,\cdot\,|w)$ is a stochastic kernel on $\Y$ given $\W,$ which determines the distribution of the observable part $y_0$ of the initial state, which may depend on the value of unobservable component $w_0=w$ of the initial state;

(v)
$c:\,\W\times\Y\times\A\mapsto  \overline{\R}_+=[0,+\infty]$ is a \textit{one-step cost function}; Dynkin and  Yushkevich~\cite[Chapter~8]{DY},
Rhenius~\cite{Rh}, Yushkevich~\cite{Yu} (see also Rieder \cite{Ri}
and B\"auerle and Rieder~\cite{BR} for a version of this model
with transition probabilities having densities)

The corresponding stochastic sequence with incomplete information \textit{evolves} as follows.
At time $t=0$ the unobservable component $w_0$ of the initial state has a given
prior distribution $p\in \P(\W).$ Let $y_0$ be the observable part
of the initial state. At
each time epoch $t=0,1,\ldots,$ if the state of the system is
$(w_t,y_t)\in\W\times\Y$ and the decision-maker chooses an action $a_t\in \A$,
then the cost $c(w_t,y_t,a_t)$ is incurred and
 the system moves to state $(w_{t+1},y_{t+1})$ according to
the transition law $P(\,\cdot\,|w_t,y_t,a_t).$

Define the \textit{observable histories}: $h_0:=y_0\in\H_0$
and
$h_t:=(y_0,a_0,y_1,a_1,\ldots,y_{t-1}, a_{t-1}, y_t)\in\H_t$ for all
$t=1,2,\dots,$
where $\H_0:=\Y$ and $\H_t:=\H_{t-1}\times \A\times \Y$
if $t=1,2,\dots$. Then a \textit{policy} for the MDPII is defined as
a sequence $\pi=\{\pi_t\}$ such that, for each $t=0,1,\dots,$
$\pi_t$ is a transition kernel on $\A$ given $\H_t$. Moreover, $\pi$
is called \textit{nonrandomized} if each probability measure
$\pi_t(\,\cdot\,|h_t)$ is concentrated at one point. The \textit{set of all policies} is denoted by $\Pi$.
The Ionescu Tulcea theorem (Bertsekas and Shreve \cite[pp.
140-141]{BS} or Hern\'andez-Lerma and Lassere
\cite[p.178]{HLerma1}) implies that a policy $\pi\in \Pi,$
initial distribution $p\in \P(\W)$
together with the transition
kernel $P$ determine a unique probability measure
$P_{p}^\pi$ on the set of all trajectories
$\mathbb{H}_{\infty}=(\W\times\Y\times \mathbb{A})^{\infty}$
endowed with the product of $\sigma$-field defined by Borel
$\sigma$-fields of $\W$, $\Y$, and $\mathbb{A}$ respectively.
The expectation with respect to this probability measure is
denoted by $\E_{p}^\pi$.

Let us specify the performance criterion. For a finite horizon
$T=0,1,\ldots,$ and for a policy $\pi\in\Pi$, let  the
\textit{expected total discounted costs} be
\begin{equation}\label{eq1}
v_{T,\alpha}^{\pi}(p):=\mathbb{E}_{p}^{\pi}\sum\limits_{t=0}^{T-1}\alpha^tc(w_t,y_t,a_t),\quad
p\in \P(\W),
\end{equation}
where $\alpha\ge 0$ is the discount factor,
$v_{0,\alpha}^{\pi}(p)=0.$ When $T=\infty$, (\ref{eq1}) defines an
\textit{infinite horizon expected total discounted cost}, and we
denote it by $v_\alpha^\pi(p).$
For any function $g^{\pi}(p)$, including
$g^{\pi}(p)=v_{T,\alpha}^{\pi}(p)$ and
$g^{\pi}(p)=v_{\alpha}^{\pi}(p)$ define the \textit{optimal value}
$g(p):=\inf\limits_{\pi\in \Pi}g^{\pi}(p),$ $p\in\P(\W).$
For a given initial distribution $p\in\P(\W)$ of the initial unobservable component $w_0,$ a policy $\pi$ is called \textit{optimal} for the respective
criterion, if $g^{\pi}(p)=g(p)$ for all $p\in \P(\W).$ A policy is called
\emph{$T$-horizon discount-optimal} if $v_{T,\alpha}^\pi(p)=v_{T,\alpha}(p),$ and it
is called \emph{discount-optimal} if $v_{\alpha}^\pi(p)=v_{\alpha}(p)$ for all $p\in \P(\W).$

We recall that a Markov decision process (MDP) is defined by its state space, action space, transition probabilities, and one-step costs. An MDP is a particular case of an MDPII. Formally speaking, an MDP $(\X,\A,P,c)$ is an MDPII $(\W\times\Y,\A,P,c)$ with $\W$ being a singleton and $\Y=\X,$ where we follow the convention that $\W\times \X=\X$ in this case. For an MDP all states are observable. If the initial state distribution $p,$ which is a probability on $\X$, is concentrated at $x\in\X,$ we write  $g^{\pi}(x)$ and $g(x)$ instead of $g^{\pi}(p)$ and $g(p).$ For an MDP, a nonrandomized policy is called \textit{Markov}, if
all decisions depend only on the current state and time. A Markov
policy is called \textit{stationary}, if all decisions depend only
on current states.

\section{Reduction of MDPIIs to MDPCIs}\label{subsec:Reduction}

In this section we formulate the well-known reduction of an MDPII $(\W\times\Y,\A,P,c)$ to an MDP called the \textit{Markov decision process with complete information (MDPCI or belief-MDP)} (\cite{BS, DY, HLerma1, Rh, Yu}). For epoch $t=0,1,\ldots$ consider
the joint conditional probability $R(B\times C|z_t,y_t,a_t)$ on next state $(w_{t+1},y_{t+1})$ given the current state $(z_t,y_t)$ and the current control action $a_t:$
\begin{equation}\label{3.3}
R(B\times C|z,y,a):=\int_{\W}P(B\times C|w,y,a)z(dw),
\end{equation}
$B\in \B(\W),$ $C\in\B(\Y),$ $(z,y,a)\in\P(\W)\times\Y\times\A.$ According to Bertsekas and Shreve \cite[Proposition~7.27]{BS}, there exists a stochastic kernel $H(z,y,a,y')[\,\cdot\,]=H(\,\cdot\,|z,y,a,y')$ on $\W$ given
$\P(\W)\times\Y\times \A\times\Y$ such that
\begin{equation}\label{3.4}
R(B\times C|z,y,a)=\int_CH(B|z,y,a,y')R(\W,dy'|z,y,a),
\end{equation}
$B\in \B(\W),$ $C\in\B(\Y),$ $(z,y,a)\in\P(\W)\times\Y\times\A.$ The stochastic kernel $H(\,\cdot\,|z,y,a,y')$ introduced in \eqref{3.4} defines a
measurable mapping $H:\,\P(\W)\times\Y\times \A\times \Y \mapsto\P(\W).$ Moreover, the mapping $y'\mapsto H(z,y,a,y')$ is defined $R(\W,\,\cdot\,|z,y,a)$-a.s. uniquely for each triple $(z,y,a)\in \P(\W)\times\Y\times\A.$

 The MDPCI is defined as an MDP
$(\P(\W)\times\Y,\A,q,\c),$ where 

(i) $\P(\W)\times\Y$ is the state space;

(ii)   $\A$ is the
action set available at all state $(z,y)\in\P(\W)\times\Y;$

(iii) $q$ on $\P(\W)\times\Y$
given $\P(\W)\times\Y\times \A$ is a stochastic kernel which determines the
distribution of the new
state as follows:
for $(z,y,a)\in \P(\W)\times\Y\times\A$
and for $D\in \B(\P(\W))$ and $C\in \B(\Y),$
\begin{equation}\label{3.7}
\begin{aligned} q(D\times C&|z,y,a):=\\ &\int_{C}\h\{H(z,y,a,y')\in D\}
R(\W,dy'|z, y, a),\end{aligned}
\end{equation}
where $\h B$ denotes the \textit{indicator of an event} $B;$

(iv) the
 one-step cost function $\c:\P(\W)\times\Y\times\A\mapsto\overline{\R}$, defined as
\begin{equation}\label{eq:c}
\c(z,y,a):=\int_{\W}c(w,y,a)z(dw),
\end{equation}
$z\in\P(\W),\,y\in\Y,\, a\in\A.$

Additional details can be found in Yushkevich \cite{Yu}, Bertsekas and Shreve \cite[Corollary~7.27.1, p.~139]{BS}, or
Dynkin and Yushkevich \cite[p.~215]{DY} . Note that the measurable particular
choice of stochastic kernel $H$ from (\ref{3.4}) does not effect 
the definition of $q$ from (\ref{3.7}).

In Section~\ref{sec:POMDPappl} we provide sufficient conditions for the existence of an optimal policy in the MDPCI $(\P(\W)\times\Y,\A,q,\c)$ in terms of the assumptions on the initial MDPII $(\W\times\Y,\A,P,c)$ and apply the results to  Platzman's model and POMDPs.  In particular, under natural conditions the existence of optimal policies, validity of optimality equations, and convergence of value iterations for MDPCIs follow from Theorem~\ref{prop:dcoe}.  If $\phi$ is a stationary discount-optimal policy for the MDPCI, then the discount optimal information policy $\pi$ for the MDPII can be defined as
$\pi(y_0,\ldots,y_t)=\phi(z_t),$ where $z_t$ is the posterior distribution of the unobservable component $w_t$ of the state $x_t$ given the observations $y_0,\ldots,y_t,$ the initial distribution $p$ of $w_0,$ and if $t>0,$ the decisions $\phi (z_s)$ for $s=0,\ldots,t-1.$  Finite-horizon discount optimal policies for MDPII can be computed similarly from Markov finite-horizon optimal policies for MDPCI.

\section{Semi-Uniform
Feller Kernels and Their Significance for MDMII}

In this section we formulate the semi-uniform Feller property for stochastic kernels. Theorem~\ref{th:mainMDPII_new} states that the transition kernel $P$ of an MDMII is semi-uniform Feller if and only if the stochastic kernel $q$ for the corresponding MDPCI satisfies this property.

Let $\S_1,\S_2,$ and $\S_3$ be Borel subsets of Polish spaces, and $\Psi$ on $\S_1\times\S_2$ given $\S_3$ be a stochastic kernel.
\begin{definition}\label{defi:unifFP}
A stochastic kernel $\Psi$ on $\S_1\times\S_2$ given $\S_3$ is semi-uniform Feller if, for each sequence $\{s_3^{(n)}\}_{n=1,2,\ldots}\subset\S_3$ that converges to $s_3$ in $\S_3$ and for each bounded continuous function $f$ on $\S_1,$
\begin{equation}\label{eq:equivWTV3}
\begin{aligned}
\lim_{n\to\infty} \sup_{B\in \B(\S_2)}& \left| \int_{\S_1} f(s_1) \Psi(ds_1,B|s_3^{(n)})\right. \\ &\left.-\int_{\S_1} f(s_1) \Psi(ds_1,B|s_3)\right|= 0.
\end{aligned}
\end{equation}
\end{definition}

We recall that the marginal measure $\Psi(ds_1,B|s_3),$ $s_3\in\S_3,$ is defined in \eqref{eq:marg_new}.
The term ``semi-uniform'' is used in Definition~\ref{defi:unifFP} because the uniform property holds in \eqref{eq:equivWTV3} only with respect to the second coordinate. If the uniform property holds with respect to both  coordinates, then the stochastic kernel $\Psi$ on $\S_1\times\S_2$ given $\S_3$ is continuous in total variation, and it is sometimes called uniformly Feller \cite{Papa}. According to Corollary~\ref{cor:Corollary 5.15.}, a semi-uniform Feller stochastic kernel is weakly continuous.

\begin{theorem}\label{th:mainMDPII_new}
Let us consider an MDMII $(\W\times\Y,\A,P,c)$ and its corresponding MDPCI $(\P(\W)\times\Y,\A,q,\c)$. The transition kernel $q$ on $\P(\W)\times\Y$
given $\P(\W)\times\Y\times \A$ for the MDPCI $(\P(\W)\times\Y,\A,q,\c)$ is semi-uniform Feller if and only if the transition kernel $P$ on $\W\times\Y$ given
$\W\times\Y\times\A$ for the MDMII $(\W\times\Y,\A,P,c)$ is semi-uniform Feller.
\end{theorem}


\section{Properties of Semi-Uniform Feller Stochastic Kernels}\label{sec:K}


Let us consider some basic definitions.
\begin{definition}\label{defi:semi}
Let $\S$ be a metric space.
 A function $f:\S\mapsto \overline{\mathbb{R}}$ is called

(i) \textit{lower semi-continuous} (l.s.c.) at a point $s\in\S$ if $\mathop{\ilim}\limits_{s'\to s}f(s')\ge f(s);$

(ii) \textit{upper semi-continuous} at $s\in\S$ if $-f$ is lower semi-continuous at $s;$

(iii)
\textit{continuous} at $s\in\S$ if $f$ is both lower and upper semi-continuous at $s;$

(iv) \textit{lower / upper semi-continuous (continuous respectively) (on $\S$)} if $f$ is lower / upper semi-continuous (continuous respectively) at each $s\in\S.$
\end{definition}

For a metric space $\S,$ let  $\F(\S),$ $\lll(\S),$ and $\C(\S)$ be the spaces of all
real-valued functions, all real-valued lower semi-continuous functions, and all real-valued continuous functions respectively defined
on the metric space $\S.$ The following definitions are taken from \cite{FKL2}.

\begin{definition}\label{defi:equi}
A set $\mathtt{F}\subset \F(\S)$ of real-valued functions on a metric space $\S$ is called
(i) \textit{lower semi-equicontinuous at a point} $s\in \S$ if $\ilim_{s'\to s}\inf_{f\in\mathtt{F}} (f(s')-f(s))\ge0;$
(ii) \textit{upper semi-equicontinuous at a point} $s\in \S$ if the set $\{-f\,:\,f\in\mathtt{F}\}$ is lower semi-equicontinuous at $s\in \S;$
(iii) \textit{equicontinuous at a point $s\in\S$}, if $\mathtt{F}$ is both lower and upper semi-equicontinuous at $s\in\S,$ that is, $\mathop{\lim}\limits_{s'\to s} \mathop{\sup}\limits_{f\in\mathtt{F}} |f(s')-f(s)|=0;$
(iv) \textit{lower / upper semi-equicontinuous (equicontinuous respectively)} (\textit{on $\S$}) if it is lower / upper semi-equicontinuous (equicontinuous respectively) at all $s \in \S;$
(v) \textit{uniformly bounded (on $\S$)}, if there exists a constant $M<+\infty $ such that $ |f(s)|\le M$ for all $s\in\S$ and  for
all $f\in \mathtt{F}.$
\end{definition}

Obviously, if a set
$\mathtt{F}\subset \F(\S)$ is lower semi-equicontinuous, then
$\mathtt{F}\subset \lll(\S).$ Moreover, if a set $\mathtt{F}\subset \F(\S)$ is equicontinuous, then $\mathtt{F}\subset \C(\S).$

\subsection{WTV-continuous Stochastic Kernels}\label{subsec:KernelsInProduct}

In this subsection we introduce WTV-continuous stochastic kernels, which are more general than semi-uniform Feller, and investigate their basic properties. In particular, Theorem~\ref{th:equivWTV} provides the equivalent definition of a  WTV-continuous stochastic kernel under the additional condition of continuity in total variation of its marginal kernel. Theorem~\ref{th:concept} establishes a necessary and sufficient condition for a stochastic kernel to be WTV-continuous and to have its marginal kernel to be continuous in total variation. This condition is Assumption~\ref{AssKern}, whose stronger version was introduced in \cite[Theorem~4.4]{Steklov}. Theorem~\ref{th:extra} describes the preservation of WTV-continuity under the integration operation.

Let $\S_1,\S_2,$ and $\S_3$ be Borel subsets of Polish spaces, and let $\Psi$ on $\S_1\times\S_2$ given $\S_3$ be a stochastic kernel.
For each set $A\in\B(\S_1)$ consider the set of functions
\begin{equation}\label{eq:familyoffunctions}
\fff^\Psi_A=\{  s_3\mapsto \Psi(A\times B |s_3):\, B\in \B(\S_2)\}
\end{equation}
mapping $\S_3$ into $[0,1].$ Consider the following type of continuity for stochastic kernels on $\S_1\times\S_2$ given $\S_3.$
\begin{definition}\label{defi:wtv}
A stochastic kernel $\Psi$ on $\S_1\times\S_2$ given $\S_3$ is called \textit{WTV-continuous}, if for each $\oo \in\tau (\S_1)$
the set of
functions $\fff^\Psi_\oo$ is lower semi-equicontinuous on $\S_3.$
\end{definition}

Similarly to Parthasarathy \cite[Theorem~II.6.1]{Part}, where the necessary and sufficient conditions for weakly convergent probability
measures were considered, the following theorem
provides several useful equivalent definitions of the WTV-continuous stochastic kernels.

\begin{theorem}\label{th:equivWTV}
For a stochastic kernel $\Psi$ on $\S_1\times\S_2$ given $\S_3,$ such that the marginal kernel $\Psi(\S_1,\,\cdot\,|\,\cdot\,)$ on $\S_2$ given $\S_3$ is continuous in total variation, the following conditions are equivalent:
(a) the stochastic kernel $\Psi$ on $\S_1\times\S_2$ given $\S_3$ is semi-uniform Feller;
(b) the stochastic kernel $\Psi$ on $\S_1\times\S_2$ given $\S_3$ is WTV-continuous;
(c) if $s_3^{(n)}$ converges to $s_3$ in $\S_3,$ then for each closed set $C$ in $\S_1$
\begin{equation}\label{eq:wtvscConv}
\lim_{n\to\infty} \sup_{B\in \B(\S_2)} \left( \Psi(C \times B|s_3^{(n)})-\Psi(C \times B|s_3)\right)= 0;
\end{equation}
(d) if $s_3^{(n)}$ converges to $s_3$ in $\S_3,$ then, for each $A\in\B(\S_1)$ such that $\Psi(\S_1,\partial A|s_3)=0,$
\begin{equation}\label{eq:equivWTV2}
\lim_{n\to\infty} \sup_{B\in \B(\S_2)} | \Psi(A \times B|s_3^{(n)})-\Psi(A \times B|s_3)|= 0;
\end{equation}
(e) if $s_3^{(n)}$ converges to $s_3$ in $\S_3,$ then, for each nonnegative bounded lower semi-continuous function $f$ on $\S_1,$
\begin{equation}\label{eq:equivWTV4}
\begin{aligned}
\ilim_{n\to\infty} \inf_{B\in \B(\S_2)} \left( \int_{\S_1}\right. &f(s_1) \Psi(ds_1,B|s_3^{(n)})\\
&\left.-\int_{\S_1} f(s_1) \Psi(ds_1,B|s_3)\right)= 0.
\end{aligned}
\end{equation}
\end{theorem}


\begin{corollary}\label{cor:Corollary 5.10} A stochastic kernel $\Psi$ on $\S_1\times\S_2$ given $\S_3$ is semi-uniform Feller if and only if the marginal kernel  $\Psi(\S_1,\,\cdot\,|\,\cdot\,)$ on $\S_2$ given $\S_3$ is continuous in total variation and at least one of conditions (b)--(e) of Theorem~\ref{th:equivWTV} holds.
\end{corollary}
%

Let us consider the following assumption. According to Example~\ref{exa:stronger}, Assumption~\ref{AssKern} is weaker than combined assumptions (i) and (ii) in \cite[Theorem~4.4]{Steklov}, where the base $\tau_b^{s_3}(\S_1)$ is the same for all $s_3\in\S_3.$

\begin{Assumption}\label{AssKern}
Let for each $s_3\in\S_3$ the topology on $\S_1$ have a countable base
$\tau_b^{s_3}(\S_1)$ such that
(i) $\S_1\in\tau_b^{s_3}(\S_1);$
 (ii)  for
each finite intersection $\oo=\cap_{i=1}^ k {\oo}_{i},$ $k=1,2,\ldots,$ of sets
$\oo_{i}\in\tau_b^{s_3}(\S_1),$ $i=1,2,\ldots,k,$
the set of
functions $\fff^\Psi_\oo,$ defined in \eqref{eq:familyoffunctions}, is equicontinuous at $s_3.$
\end{Assumption}

Note that Assumption~\ref{AssKern}(ii) holds if and only if for
each finite intersection $\oo=\cap_{i=1}^ k {\oo}_{i}$ of sets
$\oo_{i}\in\tau_b^{s_3}(\S_1),$ $i=1,2,\ldots,k,$
\begin{equation}\label{eq:equivWTV0new1}
\lim_{n\to\infty} \sup_{B\in \B(\S_2)} \left| \Psi(\oo \times B|s_3^{(n)})-\Psi(\oo \times B|s_3)\right|= 0
\end{equation}
if $s_3^{(n)}$ converges to $s_3$ in $\S_3.$

The following example demonstrates that the version of Assumption~\ref{AssKern} with the same base $\tau_b(\S_1)$ for all $s_3\in\S_3$ is essentially stronger than Assumption~\ref{AssKern}.

\begin{example}\label{exa:stronger}
{\rm
Let $\S_1=\S_3:=\mathbb{R},$ ${\rm Card\,}\S_2=1,$ and $\Psi(S_1|s_3):=\h\{s_3\in S_1\}$  $\forall S_1\in \B(\S_1)$ and $\forall s_3\in\S_3.$

Let us prove that Assumption~\ref{AssKern} holds. Indeed, for a fixed $s_3\in \mathbb{R}$ let us consider the countable base $\tau_b^{s_3}(\mathbb{R})=\{\mathbb{R}\}\cup\{(a+\sqrt{2},b+\sqrt{2})\,:\,a,b\in \mathbb{Q},\,a<b\}$ for $s_3\in\mathbb{Q},$ and
$\tau_b^{s_3}(\mathbb{R})=\{\mathbb{R}\}\cup\{(a,b)\,:\,a,b\in \mathbb{Q},\,a<b\}$ for $s_3\notin\mathbb{Q},$ where $\mathbb{Q}$ is the set of rational numbers. Note that this base satisfies the following properties :
(a) $\mathbb{R}\in \tau_b^{s_3}(\mathbb{R}),$ (b) $\oo=\cap_{i=1}^ k {\oo}_{i}\in \tau_b^{s_3}(\mathbb{R})$ for any $k=1,2,\ldots$ and
$\{\oo_{i}\}_{i=1}^{k}\subset\tau_b^{s_3}(\mathbb{R}),$ and (c) $s_3\notin\partial \oo$ for all $\oo\in \tau_b^{s_3}(\mathbb{R}).$ Statement~(a) implies that Assumption~\ref{AssKern}(i) holds.
Assumption~\ref{AssKern}(ii) holds because, according to (b) each  finite intersection $\oo=\cap_{i=1}^ k {\oo}_{i}$ of sets
$\oo_{i}\in\tau_b^{s_3}(\mathbb{R}),$ $i=1,2,\ldots,k,$ belongs to $\tau_b^{s_3}(\mathbb{R}),$ and according to (c) the function $s\mapsto \h\{s\in \oo\}$ is continuous at $s_3.$ Thus, Assumption~\ref{AssKern} holds.

Assumption~\ref{AssKern} does not hold with the same base $\tau_b(\S_1)$ for all $s_3\in\S_3$ because for any nonempty open set $\oo\in\tau(\S_1)\setminus \{\S_1\}$ there exist $s_3^*\in \partial \oo$ and a sequence $\{s_3^{(n)}\}_{n=1,2,\ldots}\subset \oo$ such that $s_3^{(n)}\to s_3^*$ in $\S_3$ as $n\to\infty,$ and, therefore, $\Psi(\oo|s_3^{(n)})=\h\{s_3^{(n)}\in \oo\}=1\not\to 0=\h\{s_3^*\in \oo\}=\Psi(\oo|s_3^*)$ as $n\to\infty,$ that is, the set of
functions $\fff^\Psi_\oo$ is not equicontinuous at $s_3^*.$ \hfill $\Box$
}
\end{example}

Theorem~\ref{th:concept} shows that Assumptions~\ref{AssKern} is a necessary and sufficient condition for semi-uniform Feller continuity.

\begin{theorem}\label{th:concept}
The stochastic kernel $\Psi$ on $\S_1\times\S_2$ given $\S_3$ is semi-uniform Feller if and only if Assumption~\ref{AssKern} holds for this kernel.
\end{theorem}

Now let $\S_4$ be a Borel subset of a Polish space, and let $\Xi$ be a stochastic kernel on $\S_1\times\S_2$ given $\S_3\times\S_4.$ Consider the stochastic kernel $\Psii$ on $\S_1\times\S_2$ given\,$\P(\S_3)\times\S_4$\,defined\,by
\begin{equation}\label{eq:extra1}
\Psii(A\times B|\mu,s_4):=\int_{\S_3}\Xi(A\times B |s_3,s_4)\mu(ds_3),
\end{equation}
$A\in\B(\S_1),\,B\in\B(\S_2),\,\mu\in\P(\S_3),\,s_4\in\S_4.$

The following theorem establishes the preservation of WTV-continuity of the integration  operation in \eqref{eq:extra1}.

\begin{theorem}\label{th:extra}
The stochastic kernel $\Psii$ on
$\S_1\times\S_2$ given $\P(\S_3)\times\S_4$ is WTV-continuous
if and only if \, $\Xi$ on $\S_1\times\S_2$ given $\S_3\times\S_4$ is WTV-continuous.
\end{theorem}

\subsection{Continuity Properties of Posterior Distributions}\label{subsec:ContBR}

In this section we establish sufficient conditions for semi-uniform Feller continuity of posterior distributions. The main result of this section is Theorem~\ref{th:CBRmain}.

Let $\S_1,\S_2,$ and $\S_3$ be Borel subsets of Polish spaces, and $\Psi$ on $\S_1\times\S_2$ given $\S_3$ be a stochastic kernel. By Bertsekas and Shreve \cite[Proposition~7.27]{BS}, there exists a stochastic kernel $\Phi$ on $\S_1$ given
$\S_2\times\S_3$ such that
\begin{equation}\label{eq:CBR1}
\Psi(A\times B|s_3)=\int_B\Phi(A|s_2,s_3)\Psi(\S_1,ds_2|s_3),
\end{equation}
$A\in \mathcal{B}(\S_1),\  B\in \mathcal{B}(\S_2),\ s_3\in\S_3.$

The stochastic kernel $\Phi(\,\cdot\,|s_2,s_3)$ on $\S_1$ given
$\S_2\times\S_3$ defines a measurable
mapping $\Phi:\,\S_2\times\S_3 \to\P(\S_1),$ where
$\Phi(s_2,s_3)(\,\cdot\,)=\Phi(\,\cdot\,|s_2,s_3).$ According to Bertsekas and Shreve \cite[Corollary~7.27.1]{BS}, for each $s_3\in
\S_3$ the mapping $\Phi(\,\cdot\,,s_3):\S_2\to\P(\S_1)$ is defined
$\Psi(\S_1,\,\cdot\,|s_3)$-almost surely uniquely in $s_2\in\S_2.$
Consider the stochastic kernel
\begin{equation}\label{eq:CBR2}
\phi(D\times B|s_3):=\int_{B}\h\{\Phi(s_2,s_3)\in D\}\Psi(\S_1,ds_2|s_3),
\end{equation}
$D\in \mathcal{B}(\P(\S_1)),\ B\in\B(\S_2),\ s_3\in\S_3.$
In models for decision making with incomplete information, $\phi$ is the transition probability between belief states, which are posterior distributions of states; \eqref{3.7}. Continuity properties of $\phi$ play the fundamental role in the studies of models with incomplete information. Theorem~\ref{th:CBRmain} characterizes such properties, and this is the reason for the title of this section.

According to Bertsekas and Shreve \cite[Corollary~7.27.1]{BS},
the particular choice of
a stochastic kernel $\Phi$ satisfying (\ref{eq:CBR1}) does not effect the
definition of $\phi$ in (\ref{eq:CBR2}) because for each $s_3\in
\S_3$ the mapping $\Phi(\,\cdot\,,s_3):\S_2\to\P(\S_1)$ is defined
$\Psi(\S_1,\,\cdot\,|s_3)$-almost surely uniquely in $s_2\in\S_2.$

Consider the following assumption.

\begin{Assumption}\label{Ass:H} There exists a stochastic
kernel $\Phi$ on $\S_1$ given $\S_2\times\S_3$ satisfying
(\ref{eq:CBR1}) such that, if a sequence
$\{s_3^{(n)}\}_{n=1,2,\ldots}\subset\S_3$ converges to
$s_3\in\S_3$ as $n\to\infty,$ then there exists a
subsequence $\{s_3^{(n_k)}\}_{k=1,2,\ldots}\subset
\{s_3^{(n)}\}_{n=1,2,\ldots}$ and a measurable subset $B$ of
$\,\S_2$ such that
\begin{equation}\label{eq:CBR3}
\begin{aligned}
&\Psi(\S_1\times B|s_3)=1\quad\mbox{and}\quad\Phi(s_2,s_3^{(n_k)})\\ &\mbox{ converges weakly to }\Phi(s_2,s_3),\quad\mbox{for all }s_2\in B.
\end{aligned}
\end{equation}
In other words, the convergence in \eqref{eq:CBR3} holds $\Psi(\S_1,ds_2|s_3)$-almost
surely.
\end{Assumption}

The following theorem, which is the main result of this section, provides necessary and sufficient conditions for WTV-continuity of a stochastic kernel $\phi$ in terms of the properties of a given stochastic kernel $\Psi.$ This theorem and the results of Subsection~\ref{subsec:KernelsInProduct} provide the necessary and sufficient conditions for the semi-uniform Feller property of the MDPCIs in terms of the conditions on the transition kernel in the initial model for decision making with incomplete information.

\begin{theorem}\label{th:CBRmain}
For a given stochastic kernel $\Psi$ on $\S_1\times\S_2$ given $\S_3,$ let the marginal kernel $\Psi(\S_1,\,\cdot\,|\,\cdot\,)$ on $\S_2$ given $\S_3$ is continuous in total variation. Then the following conditions are equivalent: (a) the stochastic kernel $\Psi$ on $\S_1\times\S_2$ given $\S_3$ is semi-uniform Feller; (b) Assumption~\ref{Ass:H} holds;
(c) if a sequence
$\{s_3^{(n)}\}_{n=1,2,\ldots}\subset\S_3$ converges to
$s_3\in\S_3$ as $n\to\infty,$ then
$\rho_{\P(\S_1)}(\Phi(s_2,s_3^{(n)}),\Phi(s_2,s_3))\to0$ in probability $\Psi(\S_1,ds_2|s_3),$
where $\rho_{\P(\S_1)}$ is an arbitrary metric that induces the topology of weak convergence of probability measures on $\S_1,$ in particular, Kantorovich-Rubinshtein metric defined in \eqref{eq:KantorRubMetr};
(d) the stochastic kernel $\phi$ on $\P(\S_1)\times \S_2$ given $\S_3$ is semi-uniform Feller;
and each of these statements implies that the stochastic kernels $\Psi$ on $\S_1\times\S_2$ given $\S_3$ and $\phi$ on $\P(\S_1)\times \S_2$ given $\S_3$ are weakly continuous.
\end{theorem}

\begin{corollary}\label{cor:Corollary 5.15.} A semi-uniform Feller stochastic kernel $\Psi$ on $\S_1\times \S_2$ given $\S_3$ is weakly continuous.
\end{corollary}

\section{Markov Decision Processes with Semi-Uniform Feller Kernels
}\label{sec:MDPwithSemi-Feller}

Let $\X_W$ and $\X_Y$ are Borel subsets of Polish spaces. In this section we consider the special class of MDPs with semi-uniform Feller transition kernels, when $\X:=\X_W\times\X_Y.$ These results are important for MDPIIs with semi-uniform Feller transition kernels from Section~\ref{sec:POMDPappl}, where $\X_W:=\P(\W)$ and $\X_Y=\Y.$

For an $\overline{\mathbb{R}}$-valued
function $f,$ defined on a nonempty subset $U$ of a metric
space $\mathbb{U},$ consider the level sets
\begin{equation}\label{def-D}
\mathcal{D}_f(\lambda;U)=\{y\in U \, : \,  f(y)\le
\lambda\},\qquad \lambda\in\R.
\end{equation}  We recall that a
function $f$ is \textit{inf-compact on $U$} if all the
level sets $\mathcal{D}_f(\lambda;U)$
are compact.

For a metric space $\U$, we denote by $\K(\U)$  the family of all
nonempty compact subsets of $\U.$
\begin{definition}{\rm( \cite[Definition~1.1]{FKN})}\label{K-inf-comp}
A function $u:\S_1\times\S_2\mapsto \overline{\mathbb{R}}$ is called $\K$-inf-compact if this function is inf-compact on $K \times\S_2$ for each $K\in \K(\S_1).$
\end{definition}

The fundamental importance of $\K$-inf-compactness is that Berge's theorem on lower semicontinuity of the minimum holds with possibly noncompact action sets;  \cite[Theorem~1.2]{FKN}. In particular, it allows us to
consider the MDPII $(\W\times\Y,\A,P,c)$ with possibly noncompact action space $\A$ and unbounded one-step cost $c$ and examine the convergence of value and policy iterations for this model in Theorem~\ref{th:mainMDPII}, in particular, for Platzman's models and partially observable MDPs in Theorem~\ref{th:mainPOMDP} and its corollaries.

\begin{definition}\label{Ass:MeasKinf}
A Borel-measurable function $u:\S_1\times\S_2\times\S_3\mapsto \overline{\mathbb{R}}$ is called measurable $\K$-inf-compact ($\M\K$-inf-compact) if for each $s_2\in \S_2$ the function $(s_1,s_3)\mapsto u(s_1,s_2;s_3)$ is $\K$-inf-compact on $\S_1\times\S_3.$
\end{definition}

As follows from Theorem~\ref{th:PresC} below, the function $\c:\P(\S_1)\times\S_2\times\S_3\mapsto\overline{\R}$ defined in \eqref{eq:c} is $\M\K$-inf-compact
if $c:\,\S_1\times\S_2\times\S_3\mapsto  \overline{\R}_+$ is an $\M\K$-inf-compact function.

Consider a discrete-time MDP $(\X,\A,q,c)$ with a state space $\X=\X_W\times \X_Y,$ an
action space $\mathbb{A},$ one-step costs $c,$ and
transition probabilities $q.$ Assume that $\X_W,\X_Y,$ and $\mathbb{A}$ are
\textit{Borel subsets} of Polish spaces. Let $LW(\mathbb{X})$ be the class of all nonnegative  Borel-measurable
functions $\varphi:\mathbb{X}\to\overline{\mathbb{R}}$
such that $w\mapsto \varphi(w,y)$ is lower semi-continuous on $\X_W$ for each $y\in\X_Y.$ For any $\alpha\ge 0$ and $u\in LW(\X),$ we consider
\begin{equation}\label{e:defeta}
\eta_u^\alpha(x,a)=c(x,a)+\alpha\int_\X  u(\tilde{x})q(d\tilde{x}| x,a),
\end{equation}
$(x,a)\in \X\times\A.$

The following theorem is the main result of this
section. It is applied in Theorem~\ref{th:mainMDPII} to MDPCIs $(\P(\W)\times\Y,\A,q,\c)$.
\begin{theorem}\label{prop:dcoe}{\rm(Expected Total Discounted  Costs)} Let us consider an MDP $(\X,\A,q,c)$ with
$\X=\X_W\times \X_Y,$ $q$ on $\X$ given $\X\times\A$ being semi-uniform Feller, and the nonnegative function $c:\X\times \A\mapsto \overline{\mathbb{R}}$ being
 $\M\K$-inf-compact.
Then
\begin{itemize}
\item[{\rm(i)}] the functions $v_{t,\alpha},$ $t=0,1,\ldots,$ and $v_\alpha$
belongs to $LW(\mathbb{X}),$ and $v_{t,\alpha}(x)\uparrow
v_\alpha (x)$ as $t \to +\infty$ for all $x\in \X;$
\item[{\rm(ii)}] $v_{t+1,\alpha}(x)=\min\limits_{a\in \A}\eta_{v_{t,\alpha}}^\alpha(x,a),$
$x\in \X,$ $t=0,1,...,$ where $v_{0,\alpha}(x)=0$ for all $x\in \X,$ and the nonempty sets $A_{t,\alpha}(x):=\{a\in \A\,:\,v_{t+1,\alpha}(x)=\eta_{v_{t,\alpha}}^\alpha(x,a) \},$ $x\in \X,$ $t=0,1,\ldots,$ satisfy the following properties: (a) the graph
${\rm Gr}(A_{t,\alpha})=\{(x,a)\,:\, x\in\X,\, a\in A_{t,\alpha}(x)\},$
$t=0,1,\ldots,$ is a Borel subset of $\X\times \mathbb{A},$ and (b)
if $v_{t+1,\alpha}(x)=+\infty,$ then $A_{t,\alpha}(x)=\A$ and, if
$v_{t+1,\alpha}(x)<+\infty,$ then $A_{t,\alpha}(x)$ is compact;
\item[{\rm(iii)}] for any $T=1,2,\ldots,$ there exists a Markov optimal
$T$-horizon policy $(\phi_0,\ldots,\phi_{T-1})$ and if, for an
$T$-horizon Markov policy $(\phi_0,\ldots,\phi_{T-1})$ the
inclusions $\phi_{T-1-t}(x)\in A_{t,\alpha}(x),$ $x\in\X,$
$t=0,\ldots,T-1,$ hold then this policy is $T$-horizon optimal;
\item[{\rm(iv)}] for $\alpha\in [0,1)$ $v_{\alpha}(x)=\min\limits_{a\in \A}\eta_{v_{\alpha}}^\alpha(x,a),$ $x\in \X,$
and the nonempty sets $A_{\alpha}(x):=\{a\in
\A\,:\,v_{\alpha}(x)=\eta_{v_{\alpha}}^\alpha(x,a) \},$ $x\in \X,$
satisfy the following properties: (a) the graph ${\rm
Gr}(A_{\alpha})=\{(x,a)\,:\, x\in\X,\, a\in A_\alpha(x)\}$  is a Borel
subset of $\X\times \mathbb{A},$ and (b) if $v_{\alpha}(x)=+\infty,$
then $A_{\alpha}(x)=\A$ and, if $v_{\alpha}(x)<+\infty,$ then
$A_{\alpha}(x)$ is compact.
\item[{\rm(v)}] for an infinite-horizon there exists a stationary
discount-optimal policy $\phi_\alpha,$ and a stationary policy is
optimal if and only if $\phi_\alpha(x)\in A_\alpha(x)$ for all
$x\in \X.$
%
\end{itemize}
\end{theorem}

\section{Total-cost Optimal Policies for MDPII and Corollaries for Platzman's  Model and Partially Observable MDPs}\label{sec:POMDPappl}

Let us state Theorem~\ref{th:mainMDPII} which is the corollary to Theorems~\ref{th:mainMDPII_new}, \ref{th:CBRmain}, and to the reduction of  MDPIIs to  MDPCIs. Then we consider Platzman's  model and POMDPs.

\begin{theorem}\label{th:mainMDPII}
Let $(\W\times\Y,\A,P,c)$ be an MDPII and $(\P(\W)\times\Y,\A,q,\c)$ be its MDPCI. Then the following conditions are equivalent :
\begin{itemize}
\item[{\rm(a)}] Assumption~\ref{AssKern} holds with $\S_1:=\W,$ $\S_2:=\Y,$ $\S_3:=\W\times\Y\times\A,$ and $\Psi:=P;$
\item[{\rm(b)}] the stochastic kernel $P$ on $\W\times\Y$ given
$\W\times\Y\times\A$ is semi-uniform Feller;
\item[{\rm(c)}] the stochastic kernel $R(\W,\,\cdot\,|\,\cdot\,)$ on $\Y$ given $\P(\W)\times\Y\times \A$ is continuous in total variation, and at least one of equivalent statements (a)--(d) of\, Theorem~\ref{th:CBRmain} with $\S_1:=\W,$ $\S_2:=\Y,$ $\S_3:=\P(\W)\times\Y\times\A,$ $\Psi:=R,$ $\Phi:=H,$ and $\phi:=q$ holds;
\item[{\rm(d)}] the stochastic kernel $q$ on $\P(\W)\times\Y$
given $\P(\W)\times\Y\times \A$ is semi-uniform Feller;
\end{itemize}
and each of them implies that the stochastic kernel $q$ on $\P(\W)\times\Y$
given $\P(\W)\times\Y\times \A$ is weakly continuous. Moreover, if nonnegative function $c$ is $\M\K$-inf-compact and one of the equivalent conditions (a)--(d) holds, then assumptions and all the conclusions of Theorem~\ref{prop:dcoe} for the MDPCI $(\P(\W)\times\Y,\A,q,\c)$ hold.
\end{theorem}

According to \cite{Rh, Yu}, for each optimal policy for the MDPCI $(\P(\W)\times\Y,\A,q,\c)$ there constructively exists an optimal policy in the original MDPII $(\W\times\Y,\A,P,c)$.
\cite[Theorem~4.4]{Steklov} establishes weak continuity of the transition kernel for the MDPCI under the more restrictive assumption than statement~(a) of Theorem~\ref{th:mainMDPII} when the countable base in this condition does not depend on the argument $s_3=(w,y,a);$ see also Example~\ref{exa:stronger}. Moreover, for any $T=1,2,\ldots$ and $\alpha\ge0,$ the value functions $\tilde{V}_{T,\alpha}(z,y),\tilde{V}_\alpha(z,y)$ in the MDPCI $(\P(\W)\times\Y,\A,q,\c),$ under the conditions of Theorem~\ref{th:mainMDPII}, are concave by $z\in \P(\W)$ for any fixed $y\in\Y,$


The proof of Theorem~\ref{th:mainMDPII} uses the following preservation property for $\M\K$-inf-compactness.

\begin{theorem}\label{th:PresC}
If $c:\,\W\times\Y\times\A\mapsto  \overline{\R}_+$ is an $\M\K$-inf-compact function, then the function $\c:\P(\W)\times\Y\times\A\mapsto\overline{\R}$ defined in \eqref{eq:c} is $\M\K$-inf-compact too.
\end{theorem}

A particular case of an MDPII is  Platzman's model introduced in \cite{Plat}.
\begin{definition}\label{defi:Plpl}
Platzman's  model is an MDPII $(\W\times\Y,\A,P,c),$ where $P$ is a stochastic kernel on $\W\times\Y$ given
$\W\times\A.$
\end{definition}

 Platzman's  model is an MDPII with the transition kernel $P(\,\cdot\,|w,y,a),$ which does not depend
on $y$.  Theorem~\ref{th:mainMDPII} implies the following corollary.

\begin{corollary}\label{cor:Plpl}
Let $(\W\times\Y,\A,P,c)$ be a Platzman's  model. Then the stochastic kernel $P$ on $\W\times\Y$ given
$\W\times\A$ is semi-uniform Feller if and only if at least one of the equivalent conditions~(a), (c), or (d) of Theorem~\ref{th:mainMDPII} holds. Moreover, if the nonnegative function $c$ is $\M\K$-inf-compact and $P$ on $\W\times\Y$ given $\W\times\A$ is semi-uniform Feller, then all assumptions and statements of Theorem~\ref{th:mainMDPII} hold.
\end{corollary}

As noticed in \cite{Plat}, the special cases of Platzman's model are two MDPs defined below, which we denote as ${\rm POMDP}_1$ and ${\rm POMDP}_2;$ see Definitions~\ref{defi:POMDP1}, \ref{defi:POMDP2}.


Let $i=1,2,$ $\X,$ $\Y,$ and $\A$ be Borel subsets of Polish spaces, $P_i(dx'|x,a)$ be a stochastic kernel on
$\X$ given $\X\times\A,$ $Q_i(dy| a,x)$ be a stochastic kernel on
$\Y$ given $\A\times\X,$ $Q_{0,i}(dy|x)$ be a stochastic kernel on
$\Y$ given $\X,$ $p$ be a probability distribution on $\X.$
\begin{definition}\label{defi:POMDP1}
${\rm POMDP}_1$ $(\X,\Y,\A,P_1,Q_1,c)$ is specified by Platzmann's model $(\X\times\Y,\A,P,c)$ with $P(B\times C| x,y,a):=P_1(B|x,a)Q_1(C|x,a),$ $B\in\B(\X),$ $C\in\B(\Y),$ $x\in\X,$ $y\in\Y,$ $a\in\A.$
\end{definition}
Let $(\X,\Y,\A,P_1,Q_1,c)$ be ${\rm POMDP}_1.$ Then, the stochastic kernel $R$ on $\X\times\Y$ given
$\P(\X)\times\Y\times \A$ specified in \eqref{3.3} takes the following form,
\begin{equation}\label{eq:3.3POMDP1}
R(B\times C|z,y,a):=\int_{\X} Q_1(C|x,a)P_1(B|x,a)z(dx),
\end{equation}
$B\in \mathcal{B}(\X),$ $C\in \mathcal{B}(\Y),$ $x\in\X$, $z\in\P(\X),$ $a\in \A.$
\begin{definition}\label{defi:POMDP2}
${\rm POMDP}_2$ $(\X,\Y,\A,P_2,Q_2,c)$ is specified by Platzmann's model $(\X\times\Y,\A,P,c)$ with \[
P(B\times C| x,y,a):=\int_B Q_2(C|a,x')P_2(dx'|x,a),
\]
$B\in\B(\X),C\in\B(\Y),\,x\in\X,\,y\in\Y,\,a\in\A.$
\end{definition}
For a ${\rm POMDP}_2$ $(\X,\Y,\A,P_2,Q_2,c)$ the stochastic kernel $R$ on $\X\times\Y$ given
$\P(\X)\times\Y\times \A$ specified in \eqref{3.3} takes the following form,
\begin{equation}\label{eq:3.3POMDP2}
R(B\times C|z,y,a):=\int_{\X}\int_B Q_2(C|a,x')P_2(dx'|x,a)z(dx),
\end{equation}
$B\in \mathcal{B}(\X),$ $C\in \mathcal{B}(\Y),$ $z\in\P(\X),$ $a\in \A.$
${\rm POMDP}_1$ is Platzman's model with observations $y_{t+1}$ being ``random functions'' of $w_{t}$ and $a_{t},$ and ${\rm POMDP}_2$ is Platzmann's  model, when observations
$y_{t+1}$ being ``random functions'' of $a_{t}$ and $w_{t+1}.$
Let us apply Theorem~\ref{th:mainMDPII} to ${\rm POMDP}_1$ and ${\rm POMDP}_2.$

In the following theorem and its corollaries we establish the semi-uniform Feller property for the
stochastic kernel $q$.

\begin{theorem}\label{th:mainPOMDP}
Let $i=1,2.$ The following conditions are equivalent for a ${\rm POMDP}_i$ :
\begin{itemize}
\item[{\rm(a)}] the stochastic kernel $P$ on $\X\times\Y$ given
$\X\times\A$ is semi-uniform Feller;
\item[{\rm(b)}] Assumption~\ref{AssKern} holds with $\S_1:=\X,$ $\S_2:=\Y,$ $\S_3:=\X\times\A,$ and $\Psi:=P;$
\item[{\rm(c)}] the stochastic kernel $R(\X,\,\cdot\,|\,\cdot\,)$ on $\Y$ given $\P(\X)\times\Y\times \A$ is continuous in total variation, and one of equivalent statements (a)--(d) of\, Theorem~\ref{th:CBRmain} holds with $\S_1:=\X,$ $\S_2:=\Y,$ $\S_3:=\P(\X)\times\Y\times\A,$ $\Psi:=R,$ $\Phi:=H,$ and $\phi:=q;$
\item[{\rm(d)}] the stochastic kernel $q$ on $\Y\times\P(\X)$ given $\X\times\A$ is semi-uniform Feller;
\end{itemize}
and each of them implies that the stochastic kernel $q$ on $\P(\X)\times\Y$
given $\P(\X)\times\Y\times \A$ is weakly continuous. Moreover, if the nonnegative function $c:\X\times\Y\times\A\mapsto\overline{\R}$ is $\M\K$-inf-compact, and one of the equivalent conditions (a)--(d) holds, then assumptions and all the conclusions of Theorem~\ref{th:mainMDPII} hold for the MDPII $(\P(\X)\times\Y,\A,q,\c)$.
\end{theorem}

\begin{corollary}\label{cor:POMDPmain1}
Let us consider ${\rm POMDP}_1.$ The stochastic kernel $P_1$ on
$\X$ given $\X\times\A$ is weakly continuous and the stochastic kernel $Q_1$ on
$\Y$ given $\X\times\A$ is continuous in total variation if and only if condition~(a) of Theorem~\ref{th:mainPOMDP} hold and, then, all the conclusions of Theorem~\ref{th:mainPOMDP} hold.
\end{corollary}

\begin{corollary}\label{cor:POMDPmain}
Let us consider ${\rm POMDP}_2.$ Each of the following assumptions :
\begin{itemize}
\item[{\rm(a)}] the stochastic kernel $P_2(\cdot|x,a)$ on $\X$ given $\X\times\A$ is
weakly continuous, and the stochastic kernel $Q_2(\cdot|a,x')$ on $\Y$
given $\A\times\X$ is continuous in
 total  variation; cp. \cite[Theorem~3.6]{FKZ};
\item[{\rm(b)}] the stochastic kernel $P_2(\cdot|x,a)$ on $\X$ given $\X\times\A$ is
continuous in total  variation, and the observation kernel $Q_2(\cdot|a,x')$ on $\Y$ given $\A\times\X$ is continuous in $a$ in total  variation; cp. \cite[Theorem~2]{Saldi};
\end{itemize}
implies condition~(a) of Theorem~\ref{th:mainPOMDP} and  all conclusions of Theorem~\ref{th:mainPOMDP}.
\end{corollary}

Corollary~\ref{cor:POMDPmain} is a more general statement than \cite[Theorem~3.6]{FKZ} and \cite[Theorem~2]{Saldi}.  The first difference is that the cost function in \cite{FKZ} does not depend on observations, and the corresponding cost function $c(x,a)$ is $\K$-inf-compact.  In this paper the cost function  $c(x,y,a)$ may depend on the observation $y.$ If $c(x,y,a)$ is $\M\K$-compact and $y$ is a dummy variable, then the function $c(x,a)=c(x,y,a)$ is $\K$-inf-compact.  The second difference is that in \cite[Theorem~2]{Saldi} the observation kernel is $Q_2(\,\cdot\,|x),$ and it  does not depend on actions. The only assumption in \cite[Theorem~2]{Saldi} is continuity of $P_2(\,\cdot\,|x,a)$ in total variation.  The conclusion of   \cite[Theorem~2]{Saldi} is weak continuity of the filter transition probability, which follows from  Theorem~\ref{th:mainPOMDP} and Corollary~\ref{cor:Corollary 5.15.}.

Let us consider Platzman's plant  $(\W\times\Y,\A,P,c)$ with the cost function $c$ that does not depend on observations $y,$ that is, $c(w,y,a)=c(w,a).$  In this case the MDPCI $(\P(\W)\times\Y,\A,q,\c)$ can be reduced to a smaller MDP $(\P(\W),\A,{\hat q},{\hat c})$ with the state space $\P(\W);$ action space $\A;$ the one-step cost  function ${\hat c}:\P(\W)\times\A\mapsto\overline{\R}$, defined as
\begin{equation}\label{eq:chat}
\hat{c}(z,a):=\int_{\W}c(w,a)z(dw),
\end{equation}
$z\in\P(\W),\, a\in\A,$ and
 the transition probability
\begin{equation}\label{3.777}
\begin{aligned} {\hat q}(D|z,a):=\int_{\Y}\h\{H(z,a,y')\in D\}
R(\W,dy'|z, a),\end{aligned}
\end{equation}
for $(z,a)\in \P(\W)\times\A$
and for $D\in \B(\P(\W)),$ where the function $H$ is defined in \eqref{3.4}, and it does not depend on $y$ because the function $R$ defined in \eqref{3.3} does not depend on $y$ for Platzman's model; Formulae \eqref{eq:chat} and \eqref{3.777} follow from
\eqref{eq:c} and  \eqref{3.7} respectively.  The described reduction of an MDPCI $(\P(\W)\times\Y,\A,q,\c)$ to the MDP $(\P(\W),\A,{\hat q},{\hat c})$ holds \cite[Theorem 2]{F2005} because in the MDPCI transition probabilities from  states $(z_t,y_t)\in (\P(\W)\times\Y$ to states $z_{t+1}\in\P(\W)$  and  costs $c(z_t,a)$ do not depend on $y_t.$ For POMDPs this smaller MDP is well-known; see, e.g., \cite{FKZ, HL}.

\section{CONCLUSIONS AND FUTURE WORKS}
\subsection{Conclusions}
The paper introduces the property of semi-uniform Feller continuity for transition probabilities and applies it to discrete-time stochastic control problems with incomplete information.    The analysis of this property provides new insights and results on the structure of optimal policies for POMDPs, and gives  sufficient conditions for the existence of optimal policies and convergence of value iterations.


\subsection{Future Works}
  The future work will investigate whether the developed procedures for statistical learning can be combined with deep learning for solving POMDPs. Another research direction is to try to relax the semi-uniform Feller condition. As the results of this paper show, this condition holds in the known situations when an optimal policy exists, even if this is not stated explicitly in the model, as this takes place with POMDPs.




\begin{thebibliography}{99}


\bibitem{Ao} Aoki, M. (1965) Optimal control of partially observable Markovian systems. \emph{J. Franklin Inst.} 280(5): 367--386.

\bibitem{As}
 {\AA}str\"om, K.J. (1965). Optimal control of Markov processes with incomplete state information. \emph{Journal of Mathematical Analysis and Applications} 10: 174–205.
%

\bibitem{BR} B\"auerle, N., Rieder, U. (2011) \textit{Markov Decision Processes with Applications to Finance,} Springer-Verlag, Berlin.


\bibitem{BS} Bertsekas, D.P.,  Shreve S.E. (1978) \textit{Stochastic Optimal Control: The Discrete-Time
Case,} Academic Press, New York 

\bibitem{Bil} Billingsley, P. (1968) \textit{Convergence of Probability Measures,}
Jonh Wiley, New York.

\bibitem{bogachev} Bogachev, V.I. (2007) \textit{Measure Theory, Volume II,}
Springer-Verlag, Berlin.
%
%
\bibitem{Dy}  Dynkin, E.B. (1965) Controlled random sequences. \emph{Theory
Probab. Appl.} 10(1): 1--14.

\bibitem{DY} Dynkin, E.B.,   Yushkevich AA (1979) \textit{{C}ontrolled
{M}arkov {P}rocesses,} Springer-Verlag, New York.

%
%
\bibitem{F2005} Feinberg, E.A. (2005) On essential information in sequential decision processes. Math. Meth. Oper. Res. 62, 399--410.

\bibitem{FK} Feinberg, E.A.,  Kasyanov, P.O. (2020) MDPs with Setwise Continuous Transition Probabilities, arXiv:2011.01325

\bibitem{FKL2} Feinberg, E.A.,  Kasyanov, P.O., Liang, Y. (2020) Fatou’s lemma in its classical form and Lebesgue's
convergence theorems for varying measures with applications to Markov decision processes. \emph{Theory Probab. Appl.} 65(2): 270--291.

\bibitem{FKNMOR}  Feinberg, E.A.,  Kasyanov, P.O., Zadoianchuk, N.V. (2012)   Average-cost Markov decision
processes with weakly continuous transition probabilities.
\emph{Math. Oper. Res.} 37(4): 591--607.

\bibitem{FKN} Feinberg, E.A., Kasyanov, P.O.,  Zadoianchuk, N.V.
(2013) Berge's theorem for noncompact image sets, \emph{J. Math.
Anal. Appl.} 397(1): 255--259.

\bibitem{TVP}  Feinberg, E.A.,  P.O. Kasyanov, P.O.,
Zadoianchuk, N.V. (2014)  \textit{Fatou's lemma for weakly converging probabilities}, Theory Probab. Appl., 58(4), 683--689.


%

\bibitem{Steklov} Feinberg, E.A., Kasyanov, P.O., Zgurovsky, M.Z. (2014) Convergence of probability measures and Markov decision models with incomplete information,
\emph{Proceedings of the Steklov Institute of Mathematics,} 287 (1), 96--117.

\bibitem{FKZ} Feinberg, E.A., Kasyanov, P.O., Zgurovsky, M.Z. (2016) Partially observable total-cost Markov decision processes with weakly continuous transition probabilities, \emph{Math. Oper. Res.}, 41(2), 656--681.


\bibitem{UFL} Feinberg, E.A.,  Kasyanov, P.O., Zgurovsky, M.Z. (2016) Uniform Fatou's lemma,  \emph{Journal of Mathematical Analysis and Applications}, 444(1), 550--567.


\bibitem{FKZ2021a} Feinberg, E.A., Kasyanov, P.O., Zgurovsky, M.Z. (2021) Semi-Uniform Feller Stochastic Kernels. arXiv:2107.02207

\bibitem{FKZ2021b} Feinberg, E.A., Kasyanov, P.O., Zgurovsky, M.Z. (2021) Markov decision processes with incomplete information and semi-uniform Feller transition probabilities. arXiv:2108.09232


\bibitem{HL}  Hern\'{a}ndez-Lerma, O. (1989) \textit{Adaptive Markov Control Processes,} Springer-Verlag, New York.

\bibitem{HLerma1} Hern\'{a}ndez-Lerma, O.,  Lassere, J.B. (1996)
\textit{Discrete-Time Markov Control Processes: Basic Optimality
Criteria,} Springer, New York.

%
%

\bibitem{Saldi} Kara, A.D., Saldi, N., Y\"uksel,S. (2019) Weak Feller property of non-linear filters. \emph{Systems $\&$ Control Letters} 134: 104512.


\bibitem{PT} Papadimitriou, C.H,, Tsitsiklis, J.N. (1987) The complexity of Markov decision processes. \emph{Math. Oper. Res.}, 12(3), 441--681.

\bibitem{Papa} Papanicolaou, G.C. (1978) Asymptotic analysis of stochastic equations. Rosenblatt M, ed.  \emph{Studies in Probability Theory}  Mathematical Association of America, Washington DC, 111--179.

\bibitem{Part}  Parthasarathy, K.R. (1967)  \textit{Probability Measures on Metric Spaces,} Academic Press, New York.

\bibitem{Plat} Platzman, L.K. (1980) Optimal infinite-horizon undiscounted control of finite probabilistic systems. \emph{SIAM Journal on Control and Optimization} 18(4): 362-380.

\bibitem{Rh}  Rhenius, D. (1974) Incomplete information in Markovian decision models. \emph{Ann. Statist.} 2(6): 1327-1334.

\bibitem{Ri} Rieder, U. (1975) Bayesian dynamic programming. \emph{Adv. Appl.
Probab.}  7(2): 330-348.


\bibitem{Rudin} Rudin, W. (1964) \textit{Principles of Mathematical Analysis,} Second edition, McGraw-Hill Inc,

%

\bibitem{ShP}
Shiryaev, A.N. (1967) Some new results in the theory of controlled
random processes. \textit{Transactions of the Fourth Prague
Conference on Information Theory, Statistical Decision Functions,
Random Processes} (Prague, 1965), pp.~131-201 (in Russian); Engl.
transl. in \textit{Select. Transl. Math. Statist. Probab.}
8(1969), 49-130.

\bibitem{Sh}
Shiryaev AN (1996) \emph{Probability,} Second edition,
Springer-Verlag, New York.

%
%

\bibitem{Yu} Yushkevich AA (1976) Reduction of a controlled Markov model with incomplete data to a problem with complete information
in the case of Borel state and control spaces. \emph{Theory
Probab. Appl.} 21(1): 153-158.
%

\end{thebibliography}
\end{document}